\documentclass[11pt,twoside]{article}
\usepackage{a4}
\usepackage{amssymb,amsmath,amsthm,latexsym}
\usepackage{amsfonts}
  \usepackage{amsfonts}
\usepackage{graphicx}

\newtheorem{theorem}{Theorem}[section]

\newtheorem{definition}[theorem]{Definition}

\newtheorem{lemma} [theorem]{Lemma}

\newtheorem{proposition}[theorem]{Proposition}
\newtheorem{remark}[theorem]{Remark}

\allowdisplaybreaks

\vspace{5cm}

\begin{document}
  
  \label{'ubf'}  
\setcounter{page}{1}                                 

\markboth {\hspace*{-9mm} \centerline{\footnotesize \sc
   The monoid structure of singular twisted virtual braids  }
                 }
                { \centerline                           {\footnotesize \sc  
        Madeti, Prabhakar \& Negi, Komal                                               } \hspace*{-9mm}              
               }

\vspace*{-2cm}

\begin{center}
{ 
       {\Large \textbf { \sc  The monoid structure of singular twisted virtual braids
                               }
       }
\\

\medskip

{\sc Madeti, Prabhakar$^{1}$ \& Negi, Komal$^{2}$ }\\
{\footnotesize $^{1,2}$Department of Mathematics, Indian Institute of Technology Ropar, Rupnagar, 140001, Punjab, India}\\
{\footnotesize e-mail: {\it prabhakar@iitrpr.ac.in$^{1}$,  komal.20maz0004@iitrpr.ac.in$^{2}$}}
}
\end{center}

\thispagestyle{empty}

\hrulefill

\begin{abstract}  
{\footnotesize  In this paper, we examine specific submonoids within the singular twisted virtual braid monoid $STVB_n$. Notably, we establish that the singular twisted virtual pure braid monoid $STVP_n$ serves as the kernel of an epimorphism from $STVB_n$ onto the symmetric group $S_n$. We identify the generators and defining relations for $STVP_n$. Additionally, we construct other epimorphisms from $STVB_n$ onto $S_n$, whose kernels are analogous to $STVP_n$, and determine their respective generators and defining relations. Furthermore, we demonstrate the embedding of the monoid $STVB_n$ into a group. Also, we provide the extension of the representation of the twisted virtual braid group to the representation of the singular twisted virtual braid monoid.
}
 \end{abstract}
 \hrulefill



\section{Introduction}

Singular twisted knot theory extends the twisted knot theory \cite{MO}, which itself is an extension of virtual knot theory. A key approach to studying knots involves examining braids; similarly, investigating singular twisted virtual braids aids in understanding singular twisted links. Alexander's theorem and Markov's theorem for singular twisted links and singular twisted virtual braids have been established in the paper~\cite{KM,KPS}, illustrating a coherent connection between these structures. Additionally, it has been observed that the singular twisted virtual braid on \(n\) strands forms a monoid. Both a monoid presentation and a reduced monoid presentation for the singular twisted virtual braid monoid are provided by them.

V. V. Vershinin~\cite{VV} provided an analogue of the braid group $B_n$ presentation for the singular braid monoid structure. R. Fenn, E. Keyman, and C. Rourke~\cite{REC} later proved that the singular braid monoid can be embedded in a group, which has a geometric interpretation involving singular braids with two types of singularities that cancel each other. 

In 2020, V. Bardakov et al.~\cite{BK} studied the submonoids of $SB_n$, focusing on the singular pure braid group $SP_n$ for $n = 2, 3$. They determined the generators, defined relations, and explored the algebraic structure of these groups. Moreover, they showed that the center $Z(SP_3)$ is a direct factor in $SP_3$. 
Subsequently, in a later work~\cite{BK-1}, they identified a finite set of generators and defining relations for the singular pure braid group $SP_n$, for $n \geq 3$, which is a subgroup of the singular braid group $SG_n$. Additionally, they introduced subgroups of camomile type and established that for $n \geq 5$, the singular pure braid group $SP_n$ is a camomile-type subgroup of $SG_n$.

Motivated by all these work on $SB_n$, C. Carmen et al.~\cite{CS} explored the algebraic structures of the virtual singular braid monoid, \(VSB_n\), and the virtual singular pure braid monoid, \(VSP_n\). The monoid $V SB_n$ is  the splittable extension of $V SP_n$ by the symmetric group $S_n$. Further they also constructed a representation for $VSB_n$.

Recently, Bardakov et al.~\cite{VNT} proved that for a given representation $\varphi: B_n \to G_n$ of the braid group $B_n$ for $n \geq 2$ into a group $G_n$. This representation can be extended to \( \Phi: SM_n \to A_n \), where \( SM_n \) represents the singular braid monoid and \( A_n \) is an associative algebra whose group of units includes \( G_n \). The feasibility of extending \( \Phi: SM_n \to A_n \) further to a representation \( \Phi': SB_n \to A_n \) for the singular braid group \( SB_n \) has also been investigated.

More recently, V.G. Bardakov et al.~\cite{VTKM} investigated the structure of the twisted virtual braid group, \(TVB_n\), and offered presentations for its subgroups. Building on their work, we examine various submonoids of \(STVB_n\) and provide presentations for them. Also, we address the solution for the extension problem for the monoid $STVB_n$.

This article is organized as follows:  
In Section 2, we present the preliminaries on singular twisted virtual braids necessary for this research.  
In Section 3, we outline certain submonoids of the singular twisted virtual braid monoid.  
In Section 4, we provide results related to the embedding of the singular twisted virtual braid monoid into a group and the extension of the representation of the twisted virtual braid group to the representation of the singular twisted virtual braid monoid.
\section{Preliminaries}

In this section, we present the established results concerning singular twisted virtual braids, focusing on the monoid structure exhibited by singular twisted virtual braids on \(n\) strands.

\begin{figure}[ht]
  \centering
    \includegraphics[width=2cm,height=3cm]{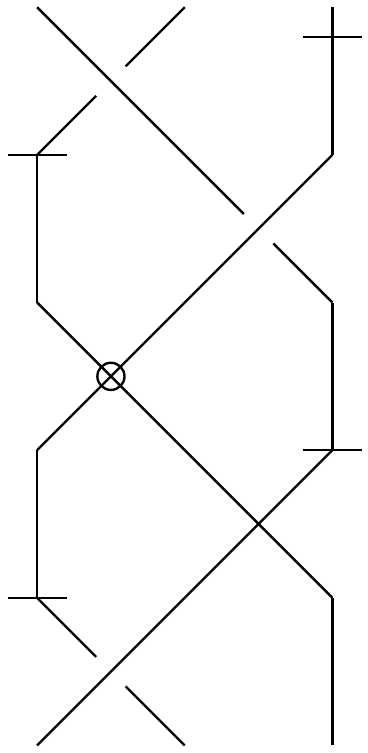}
        \caption{A singular twisted virtual braid on 3 strands}
        \label{ext}
        \end{figure}

\begin{definition}
A singular twisted virtual braid diagram with \(n\) strands (or of degree \(n\)) is a collection of \(n\) smooth or polygonal curves, referred to as strands, in \(\mathbb{R}^2\). These strands connect points \((i, 1)\) to points \((p_i, 0)\) for \(i = 1, \dots, n\), where \((p_1, \dots, p_n)\) is a permutation of \((1, \dots, n)\). The strands are required to be monotonic with respect to the second coordinate, and any intersections between the strands must be transverse double points, each labeled with information indicating whether it is a positive, negative, virtual, or singular crossing. Additionally, the strands may include bars, which are short arcs that intersect the strands transversely. For an illustration, see Figure~\ref{ext}, where the five crossings are, from top to bottom, negative, positive, virtual, singular, and positive.
\end{definition}

\begin{figure}[ht]
  \centering
    \includegraphics[width=8cm,height=3.5cm]{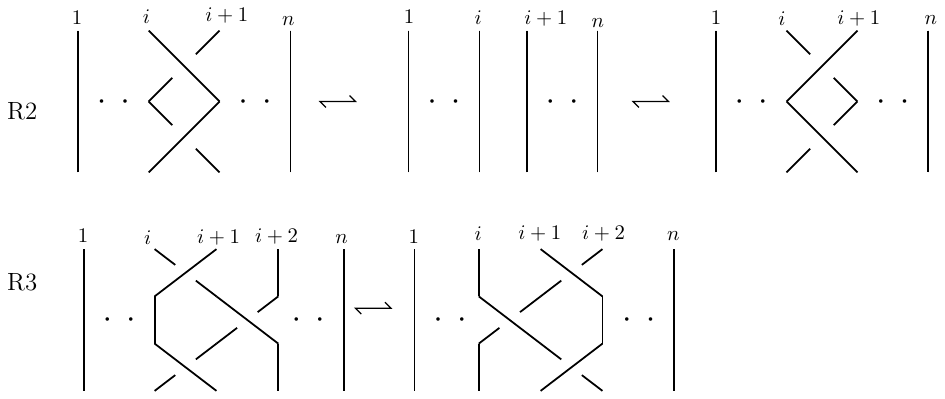}
        \caption{Classical braid moves}
        \label{bmoves}
        \end{figure}

\begin{figure}[ht]
  \centering
    \includegraphics[width=10cm,height=4cm]{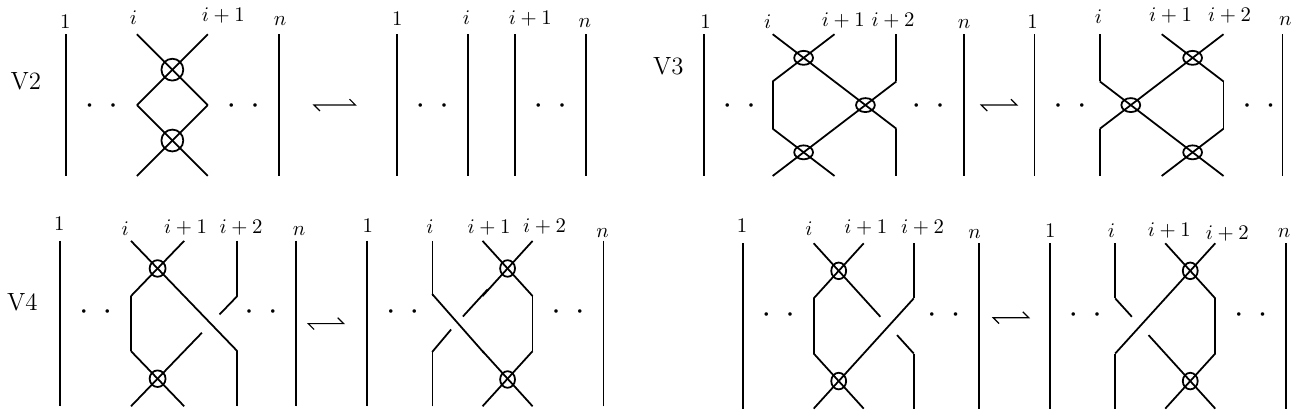}
        \caption{Virtual braid moves}
        \label{vbmoves}
        \end{figure}

\begin{figure}[ht]
  \centering
    \includegraphics[width=9cm,height=4cm]{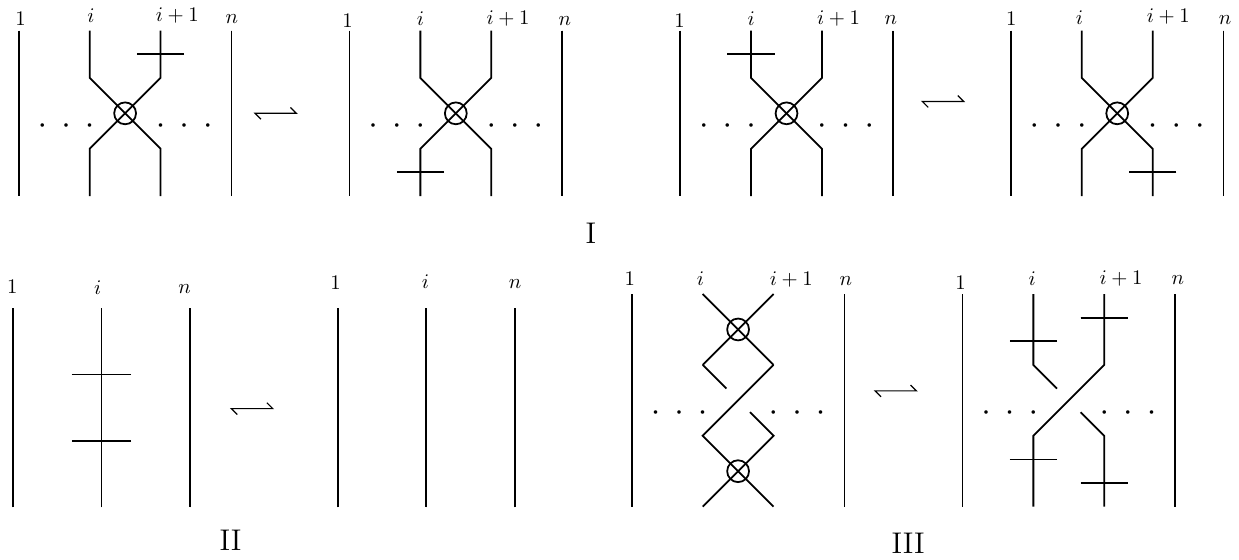}
        \caption{Twisted braid moves}
        \label{moves}
        \end{figure}

\begin{figure}[ht]
  \centering
    \includegraphics[width=9cm,height=4.5cm]{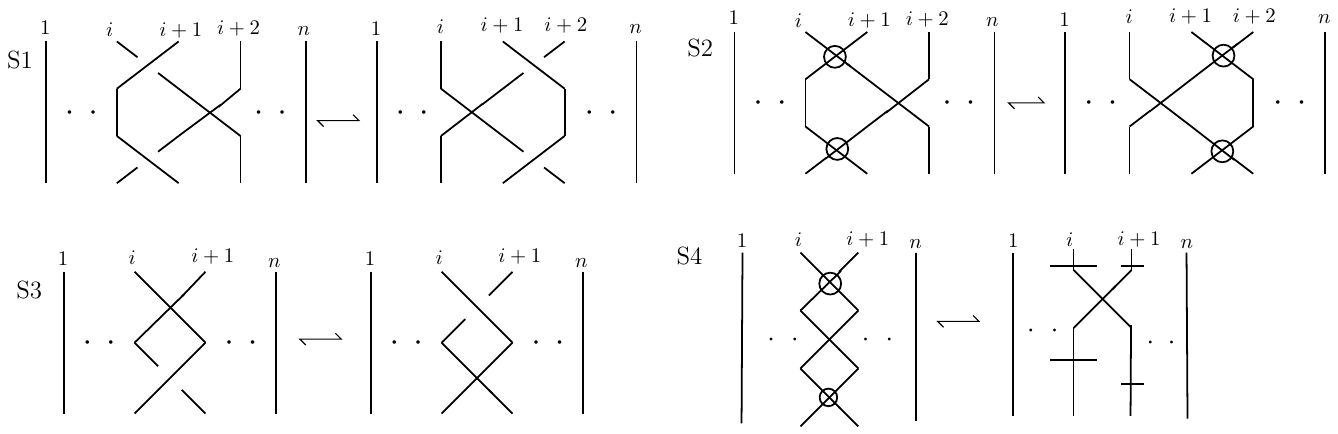}
        \caption{Singular braid moves}
        \label{smoves}
        \end{figure}

\begin{definition}
Two singular twisted virtual braid diagrams, \( b \) and \( b' \), of degree \( n \) are said to be {\it equivalent} if one can be transformed into the other using the classical, virtual, twisted, and singular braid moves illustrated in Figures~\ref{bmoves}, \ref{vbmoves}, \ref{moves}, and \ref{smoves}, respectively.
\end{definition}

\begin{definition}
   A {\it singular twisted virtual braid} refers to the equivalence class comprising singular twisted virtual braid diagrams.
\end{definition}
\subsection{Singular twisted virtual braid monoid}
In a manner similar to singular braids and virtual singular braids, the collection of singular twisted virtual braids forms a monoid, where the operation is defined by concatenation, akin to the structure of the braid group. The set of singular twisted virtual braids with \(n\) strands is also a monoid, represented by \(STVB_n\).
\begin{figure}[ht]
  \centering
    \includegraphics[width=12cm,height=2cm]{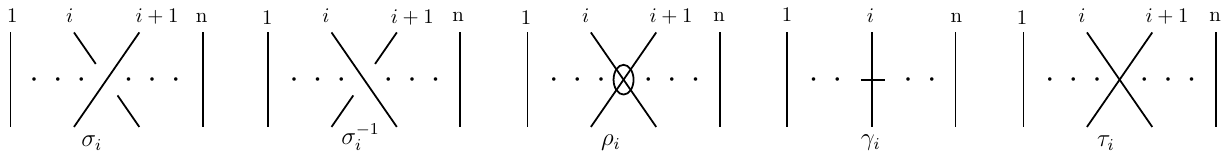}
        \caption{Generators for the monoid of singular twisted virtual braids}
        \label{gen}
        \end{figure}
\begin{theorem}~\cite{KM}\label{thm:StandardPresentation2}
The monoid $STVB_n$ is generated by standard generators, 
$\sigma_i^{\pm 1}$, $\rho_i$, $\tau_i$ $(i=1, \dots, n-1)$, and $\gamma_i$ $(i=1, \dots, n)$ satisfy
the following relations, where $e$ represents the identity element:
    \begin{align}
       \sigma_i \sigma_j & = \sigma_j \sigma_i  & \text{ for } & |i-j| > 1; \label{rel-height-ss}\\
        \sigma_i \sigma_{i+1} \sigma_i & = \sigma_{i+1} \sigma_i \sigma_{i+1} & \text{ for } & i=1,\ldots, n-2; \label{rel-sss}\\
        \rho_i^2 & = e  & \text{ for } & i=1,\ldots, n-1; \label{rel-inverse-v}\\
        \rho_i \rho_{j} & = \rho_{j} \rho_i & \text{ for } & |i-j| > 1 ; \label{rel-height-vv}\\
        \rho_i \rho_{i+1} \rho_i & = \rho_{i+1} \rho_i \rho_{i+1} & \text{ for } & i=1,\ldots, n-2; \label{rel-vvv}\\
        \sigma_i \rho_{j} & = \rho_{j} \sigma_i &  \text{ for } & |i-j| >1  ; \label{rel-height-sv}\\
        \rho_i \sigma_{i+1} \rho_i & = \rho_{i+1} \sigma_i \rho_{i+1} & \text{ for } & i=1,\ldots, n-2; \label{rel-vsv}\\
        \gamma_i^2 & = e & \text{ for } & i=1,\ldots, n; \label{rel-inverse-b}\\  
        \gamma_i \gamma_j & = \gamma_j \gamma_i   & \text{ for } & i,j=1,\ldots, n; \label{rel-height-bb} \\
        \gamma_j \rho_i & = \rho_i \gamma_j & \text{ for } & j\neq i, i+1; \label{rel-height-bv}\\
        \sigma_i\gamma_j & = \gamma_j\sigma_i & \text{ for } & j\neq i, i+1; \label{rel-height-sb}\\
        \gamma_{i+1} \rho_i & = \rho_{i} \gamma_i & \text{ for } & i=1,\ldots, n-1; \label{rel-bv} \\
        \rho_{i} \sigma_i \rho_{i} & = \gamma_{i+1} \gamma_i \sigma_{i} \gamma_i \gamma_{i+1} & \text{ for } &  i=1,\ldots, n-1; \label{rel-twist-III}\\
        \sigma_i \sigma_i^{-1} & = e  & \text{ for } & i=1,\ldots, n-1; \label{rel-height-ss1}\\
        \tau_i \tau_j & = \tau_j \tau_i & \text{ for } & |i-j| > 1 ; \label{rel-height-vv2}\\
        \sigma_i \tau_j & = \tau_j \sigma_i &  \text{ for } & |i-j| >1  ; \label{rel-height-sv1}\\
        \sigma_i \tau_{i} & = \tau_i \sigma_i &  \text{ for } &  i=1,\ldots, n; \label{rel-height-sv2}\\
        \sigma_i \sigma_{i+1} \tau_i & = \tau_{i+1} \sigma_i \sigma_{i+1} & \text{ for } & i=1,\ldots, n-2; \label{rel-sss1}\\
        \sigma_{i+1} \sigma_{i} \tau_{i+1} & = \tau_{i} \sigma_{i+1} \sigma_{i} & \text{ for } & i=1,\ldots, n-2; \label{rel-sss2}\\
        \tau_i \rho_{j} & = \rho_{j} \tau_i & \text{ for } & |i-j| > 1 ; \label{rel-height-vv1}\\
        \rho_i \tau_{i+1} \rho_i & = \rho_{i+1} \tau_i \rho_{i+1} & \text{ for } & i=1,\ldots, n-2; \label{rel-vvv1}\\
        \tau_i\gamma_j & = \gamma_j\tau_i & \text{ for } & j\neq i, i+1; \label{1rel-height-sb}\\
       \rho_{i} \tau_i \rho_{i} & = \gamma_{i+1} \gamma_i \tau_{i} \gamma_i \gamma_{i+1} & \text{ for } &  i=1,\ldots, n-1. \label{1rel-twist-III}
    \end{align}
\end{theorem}
\begin{remark}
    By taking the generators \( \sigma_i \), \( \rho_i \) (\( i = 1, \dots, n-1 \)), and \( \gamma_i \) (\( i = 1, \dots, n \)) along with the relations specified in (\ref{rel-height-ss}) to (\ref{rel-twist-III}), we obtain a group known as the twisted virtual braid group \( TVB_n \)~\cite{KPS}.
\end{remark}
Also the presentation of the monoid $STVB_n$ given in Theorem~\ref{thm:StandardPresentation2} can be reduced to a presentation with less generators and less relations by rewriting $\sigma^{\pm 1}_i$ $(i=2,\ldots, n-1)$, $\tau_i$ $(i=2,\ldots, n-1)$, and $\gamma_i$ $(i=2,\ldots, n)$ in terms of $\sigma_1^{\pm 1}$, $\tau_1$, $\gamma_1$ and $\rho_1, \dots, \rho_{n-1}$  as follows:
\begin{align}
    \sigma_i^{\pm 1} & =(\rho_{i-1}\ldots \rho_1)(\rho_i \ldots \rho_2)\sigma_1^{\pm 1}(\rho_2  \ldots \rho_i)(\rho_1  \ldots \rho_{i-1}) & \text{ for } & i=2,\ldots, n-1,  \label{1st reduction} \\ 
     \tau_i & =(\rho_{i-1}\ldots \rho_1)(\rho_i \ldots \rho_2)\tau_1(\rho_2  \ldots \rho_i)(\rho_1  \ldots \rho_{i-1}) & \text{ for } & i=2,\ldots, n-1,  \label{2nd reduction} \\
    \gamma_i & =(\rho_{i-1}\ldots \rho_1)\gamma_1(\rho_1  \ldots \rho_{i-1}) & \text{ for } & i=2,\ldots, n.  \label{3rd reduction}
\end{align}
\begin{theorem}~\cite{KM}\label{thm:ReducedPresentation}
The monoid $STVB_n$ can be described through a presentation involving generators are $\sigma^{\pm 1}_1,\tau_1, \gamma_1$,  $\rho_1,\dots, \rho_{n-1}$   
and the defining relations are (\ref{rel-inverse-v})-(\ref{rel-vvv}), and as follows:  
\begin{align}
\sigma_1(\rho_2\rho_3\rho_1\rho_2\sigma_1\rho_2\rho_1\rho_3\rho_2) & = (\rho_2\rho_3\rho_1\rho_2\sigma_1\rho_2\rho_1\rho_3\rho_2)\sigma_1, &  & \label{relB-height-ss} \\
  (\rho_1\sigma_1\rho_1)(\rho_2\sigma_{1}\rho_2)(\rho_1\sigma_1\rho_1) & = (\rho_2\sigma_1\rho_2)(\rho_1\sigma_{1}\rho_1)(\rho_2\sigma_1\rho_2), & & \label{relB-sss} \\
 \sigma_1\rho_{j}  & = \rho_{j}\sigma_1 & \text{ for } & j = 3, \ldots, n-1; \label{relB-height-sv}\\
 \gamma_1^2  & = e,  & &  \label{relB-inverse-b} \\
 \gamma_1\rho_{j} & =  \rho_{j}\gamma_1 & \text{ for } & j = 2, \ldots, n-1; & \label{relB-height-bv}\\
  \gamma_1\rho_1\gamma_1\rho_1 & =\rho_1\gamma_1\rho_1\gamma_1, & &  \label{relB-height-bb}\\
  \gamma_1\rho_1\rho_2\sigma_1\rho_2\rho_1 & = \rho_1\rho_2\sigma_1\rho_2\rho_1\gamma_1, & & \label{relB-height-sb}\\
  \gamma_{1}\rho_1\gamma_1\sigma_{1} \gamma_1\rho_1\gamma_{1} & = \sigma_1. & & \label{relB-bv}\\
  \sigma_1\sigma^{-1}_1 & =e &  \label{sigma}\\
  \sigma_1\tau_1 &=\tau_1\sigma_1 & \label{tausigma}\\
  \tau_1 \rho_i & =\rho_i \tau_1 & \text{ for }  & i=3, \ldots n-1; \label{vtau}\\
  \tau_1(\rho_1\rho_2\sigma_1\rho_2\rho_1)\sigma_1 & = (\rho_1\rho_2\sigma_1\rho_2\rho_1)\sigma_1(\rho_1\rho_2\tau_1\rho_2\rho_1) \label{1}\\
  \tau_1(\rho_2\rho_3\rho_1\rho_2\sigma_1\rho_2\rho_1\rho_3\rho_2) &=(\rho_2\rho_3\rho_1\rho_2\sigma_1\rho_2\rho_1\rho_3\rho_2)\tau_1 \label{0}\\
  \tau_1(\rho_2\rho_3\rho_1\rho_2\tau_1\rho_2\rho_1\rho_3\rho_2) &=(\rho_2\rho_3\rho_1\rho_2\tau_1\rho_2\rho_1\rho_3\rho_2)\tau_1 \label{2}\\
   \gamma_1\rho_1\rho_2\tau_1\rho_2\rho_1 & = \rho_1\rho_2\tau_1\rho_2\rho_1\gamma_1, & & \label{relB-height-sb2}\\
\gamma_{1}\rho_1\gamma_1\tau_{1} \gamma_1\rho_1\gamma_{1} & = \tau_1. & & \label{relB-bv2} 
\end{align}
  
\end{theorem}
Now, the question is can we define different submonoids of the monoid $STVB_n$ and find a presentation for them?

One possible approach is outlined in the following section.
\section{Some submonoids of $STVB_n$}  

In this section, we define various epimorphisms from \( STVB_n \) onto \( S_n \). By determining the kernel of each epimorphism, we identify the corresponding normal submonoids of \( STVB_n \).

One can define some epimorphism of $STVB_n$ onto $S_n$. The first one is
$$\varphi_1 \colon STVB_n \to S_n,~~\sigma_i^{\pm 1} \mapsto \rho_i,~~\tau_i \mapsto \rho_i, ~~\rho_i \mapsto \rho_i,~~i=1, 2, \ldots,n-1,~~\gamma_j \mapsto e,~~j=1, 2, \ldots, n.$$
Its kernel $\ker(\varphi_1)$ is the singular twisted virtual pure braid monoid $STVP_n$.

The second epimorphism  is
$$\varphi_2 \colon STVB_n \to S_n,~~\sigma_i^{\pm 1} \mapsto e,~~\tau_i \mapsto e, ~~\rho_i \mapsto \rho_i,~~i=1, 2, \ldots,n-1,~~\gamma_j \mapsto e,~~j=1, 2, \ldots, n.$$
Its kernel $\ker(\varphi_2)$ is denoted by $STVH_n$.

The third epimorphism  is
$$\varphi_3 \colon STVB_n \to S_n,~~\sigma_i^{\pm 1} \mapsto \rho_i,~~\tau_i \mapsto e, ~~\rho_i \mapsto \rho_i,~~i=1, 2, \ldots,n-1,~~\gamma_j \mapsto e,~~j=1, 2, \ldots, n.$$
Its kernel $\ker(\varphi_3)$ is denoted by $M_n$.

It is clear that $STVP_n$, $STVH_n$, and $M_n$ are normal submonoids of index $n!$ of $STVB_n$and the following short exact sequences hold:
$$1 \to STVP_n \to STVB_n \to S_n \to 1,$$
$$1 \to  STVH_n \to  STVB_n \to S_n \to 1,$$
$$1 \to M_n \to STVB_n \to S_n \to 1.$$
\subsection{Structure of singular twisted virtual pure braid monoid $STVP_n$.}
\begin{figure}[ht]
  \centering
    \includegraphics[width=8cm]{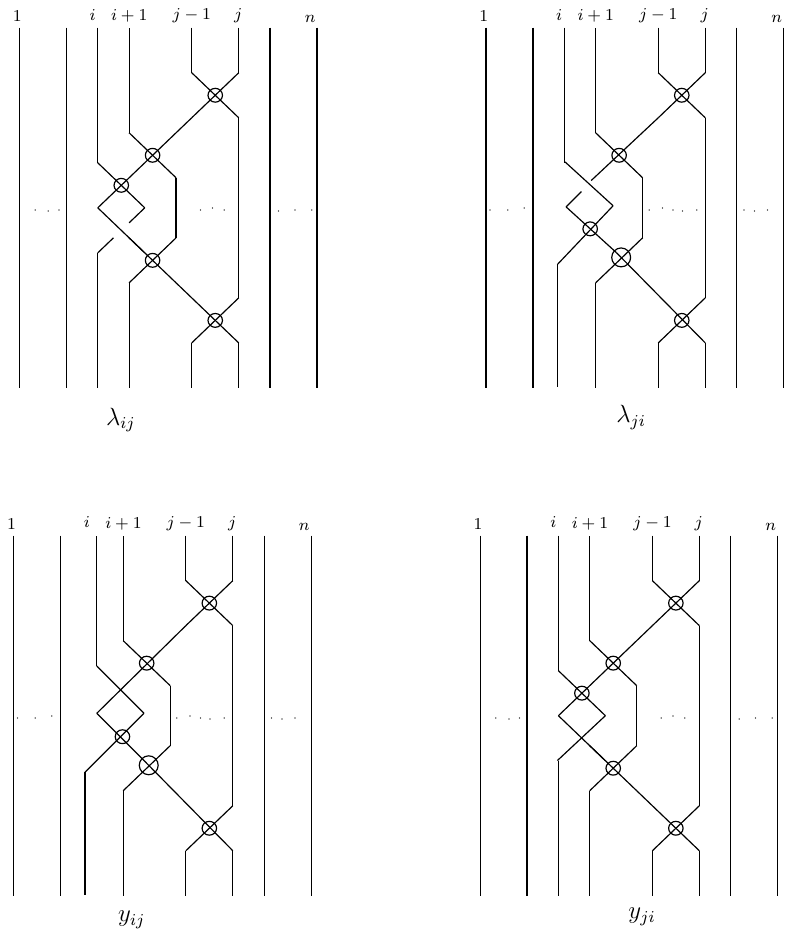}
        \caption{Generators of the monoid of singular twisted virtual pure braids}
        \label{genp}
        \end{figure}

Let us consider the elements as follows:
$$\lambda_{i,i+1}=\rho_i\sigma_i^{-1},\text{   } \lambda_{i+1,i}=\rho_i\lambda_{i,i+1}\rho_i, \text{   }  i=1, 2, \ldots, n-1,$$
$$\lambda_{i,j}=\rho_{j-1}\rho_{j-2}\cdots \rho_{i+1}\lambda_{i,i+1}\rho_{i+1}\cdots \rho_{j-2}\rho_{j-1},$$
$$\lambda_{j,i}=\rho_{j-1}\rho_{j-2}\cdots \rho_{i+1}\lambda_{i+1,i}\rho_{i+1}\cdots \rho_{j-2}\rho_{j-1},\text{   }  1 \leq i \leq j-1 \leq n-1,$$
$$ y_{i,i+1}=\tau_i\rho_i,\text{   } y_{i+1,i}=\rho_i y_{i,i+1}\rho_i, \text{   }  i=1, 2, \ldots, n-1,$$
$$y_{i,j}=\rho_{j-1}\rho_{j-2}\cdots \rho_{i+1}y_{i,i+1}\rho_{i+1}\cdots \rho_{j-2}\rho_{j-1},$$
$$y_{j,i}=\rho_{j-1}\rho_{j-2}\cdots \rho_{i+1}y_{i+1,i}\rho_{i+1}\cdots \rho_{j-2}\rho_{j-1},\text{   }  1 \leq i \leq j-1 \leq n-1,$$
which generate singular virtual pure braid monoid $SVP_n$~\cite{CS}. Geometrically, these generators are shown in Figure~\ref{genp}.

Before proving presentation of $STVP_n$, we require the following results. 

\begin{lemma}\label{gamma}
The group \( S_n \) acts on the set \( \{\gamma_i \mid 1 \leq i \leq n\} \) via conjugation, and this action is transitive. Let \( a \) denote an element of \( \langle \rho_1, \rho_2, \ldots, \rho_{n-1} \rangle \), and let \( \bar{a} \) represent its image in \( S_n \) under the isomorphism defined by \( \rho_i \mapsto (i, i+1) \) for \( i = 1, 2, \ldots, n-1 \). Then, for any generator \( \lambda_{i,j} \), \( y_{i,j} \), or \( \gamma_i \) of \( STVP_n \), the following holds.
$$a^{-1} \lambda_{i,j} a = \lambda_{(i)\bar{a}, (j)\bar{a}},~~a^{-1} y_{i,j} a = y_{(i)\bar{a}, (j)\bar{a}},\text{ and } a^{-1} \gamma_{i} a = \gamma_{(i)\bar{a}},$$
where $(k)\bar{a}$ is the image of $k$ under the action of the permutation $\bar{a}$.
\end{lemma}

To determine the generators and defining relations of \( STVP_n \), the Reidemeister-Schreier method is employed (refer to \cite[Chapter 2.3]{WAD} for details). As a Schreier set of coset representation of $STVP_n$ in $STVB_n$ we take the same set $\Lambda_n$, which is used  in finding the presentation of $TVP_n$,
$$
\Lambda_n = \left\{ \prod\limits_{k=2}^n m_{k,j_k}~ |~ 1 \leq j_k
\leq k \right\}
$$
where $m_{kl}=\rho_{k-1}\rho_{k-2}\cdots \rho_l$ for $l<k$ and $m_{kl}=1$ in the other cases.

\begin{theorem}\label{sTVP_n}
    The monoid $STVP_n$ has a presentation with the generators $\lambda_{kl}^{\pm 1},y_{kl}$, 
    $1 \leq k\neq l \leq n$, and $\gamma_j$, $1\leq j\leq n$. The defining relations are (\ref{rel-inverse-b})-(\ref{rel-height-bb}) and as follows:
     \begin{align}
     \lambda_{ij}\lambda_{ij}^{-1}&=e\label{inverseexists}\\
\lambda_{ij}\lambda_{kl} &=\lambda_{kl}\lambda_{ij},\label{comm-clas}\\
\lambda_{ki}(\lambda_{kj}\lambda_{ij}) &=(\lambda_{ij}\lambda_{kj})\lambda_{ki},\label{classical}\\
\lambda_{ij}\gamma_k & =\gamma_k\lambda_{ij},\label{g3}\\
\lambda_{ij} &=\gamma_i\gamma_j\lambda_{ji}\gamma_j\gamma_i,\label{g4}\\
y_{ij} y_{kl} &=y_{kl}y_{ij},\label{comm-clas1}\\
\lambda_{ij}y_{kl} &=y_{kl}\lambda_{ij},\label{comm-clas2}\\
\lambda_{ki}(\lambda_{kj}y_{ij}) &=(y_{ij}\lambda_{kj})\lambda_{ki},\label{classical2}\\
\lambda_{lk}y_{kl} &=\lambda_{kl}y_{lk},\label{comm-clas3}\\
y_{ki}(\lambda_{kj}\lambda_{ij}) &=(\lambda_{ij}\lambda_{kj})y_{ki},\label{classical3}\\
y_{ij}\gamma_k & =\gamma_k y_{ij},\label{g5}\\
y_{ij} &=\gamma_i\gamma_j y_{ji}\gamma_j\gamma_i,\label{g6}
\end{align}
Where each distinct letter represents a unique index.
\end{theorem}
\begin{proof}
    Define the map $\Bar{}: STVB_n \to \Lambda_n$ that assigns to each element $w \in STVB_n$ to its representative $\overline{w}$ from $\Lambda_n$. Consequently, the element $w(\overline{w})^{-1}$ belongs to $STVP_n$. According to the Theorem 2.7 of~\cite{WAD} the submonoid $STVP_n$ is generated by
    $$s_{\lambda,a}=\lambda a \cdot (\overline{\lambda a})^{-1},$$
    where $\lambda \in \Lambda_n$ and $a \in \{\rho_1, \ldots, \rho_{n-1}, \sigma^{\pm 1}_1, \ldots, \sigma^{\pm 1}_{n-1}, \tau_1,\dots, \tau_{n-1}, \gamma_1, \ldots, \gamma_ n\}$.
    
    As calculated in the proof of Theorem 17~\cite{CS} that 
    $ s_{\lambda, \rho_i} =e~~\forall \lambda, \rho_i$,
    \begin{align*}
        s_{\lambda, \sigma_i}&=\lambda (s_{e,\sigma_i} )\lambda^{-1}&=\lambda (\sigma_i\rho_i)\lambda^{-1} &=\lambda (\lambda_{i,i+1}^{-1})\lambda^{-1} & \text { (which is equal to some } \lambda_{kl}^{\pm 1}), \\
        s_{\lambda, \tau_i} &=\lambda (s_{e,\tau_i} )\lambda^{-1} &=\lambda (\tau_i\rho_i)\lambda^{-1} &=\lambda (y_{i,i+1})\lambda^{-1} & \text { (which is equal to some } y_{kl}),\\
        s_{\lambda, \sigma_i^{-1}}&=\lambda (s_{e,\sigma_i^{-1}} )\lambda^{-1}&=\lambda (\sigma_i^{-1}\rho_i)\lambda^{-1} &=\lambda (\lambda_{i+1,i})\lambda^{-1} & \text { (which is equal to some } \lambda_{kl}^{\pm 1}).
    \end{align*}
    Now, consider the new generators $$s_{\lambda,\gamma_i}=\lambda(s_{e,\gamma_i})\lambda^{-1}$$
    we know, $s_{e,\gamma_i}=\gamma_i$. So, $s_{\lambda,\gamma_i}=\lambda(\gamma_i)\lambda^{-1}$ which is equal to some $\gamma_j$ by Lemma~\ref{gamma}.
   Therefore, generators of the monoid $STVP_n$ are $\lambda_{kl}$, $\lambda_{kl}^{-1}$, $y_{kl}$ for $1 \leq k\neq l \leq n$ and, $\gamma_j$ for $1\leq j\leq n$.

  To determine the defining relations of \( STVP_n \), we introduce a rewriting process \( f \). This process allows us to rewrite a word that represents an element of \( STVP_n \) (written in terms of the generators of \( STVB_n \)) as a word in the generators of \( STVP_n \). Let us associate to reduced word $$u=a_1^{\epsilon_1}a_2^{\epsilon_2}\ldots a_v^{\epsilon_v}, \epsilon_l=\pm 1,$$ 
   where
   $a_l \in \{\sigma^{\pm 1}_1,\sigma^{\pm 1}_2, \ldots, \sigma^{\pm 1}_{n-1}, \rho_1,\rho_2, \ldots, \rho_{n-1}, \tau_1,\tau_2, \ldots, \tau_{n-1}, \gamma_1,\gamma_2, \ldots, \gamma_n\}$.
   
   The word
   $$f(u)=s_{k_1,a_1}^{\epsilon_1}s_{k_2,a_2}^{\epsilon_2}\ldots s_{k_v,a_v}^{\epsilon_v}$$
   in the generators of $STVP_n$, where $k_j$ represents the $(j-1)^{th}$ initial segment of the word $u$ if $\epsilon_j=1$, and a representative of the $j^{th}$ initial segment of $u$ if $\epsilon_j=-1$.

   According to the Theorem 2.9 in~\cite{WAD}, the monoid $STVP_n$ is defined by the relations
   $$r_{\mu, \lambda}=f(\lambda r_\mu \lambda^{-1})=\lambda f( r_\mu) \lambda^{-1}, \lambda \in \Lambda_n,$$
   where $r_\mu$ is the defining relations of $STVB_n$.

   The relations (\ref{rel-height-ss}) to (\ref{rel-vvv1}) give rise to relations (\ref{inverseexists}) and (\ref{g4}) proved in~\cite{VTKM,CS}.

   Denote $r_1:=\tau_i\gamma_j =\gamma_j\tau_i$, where $j\neq i, i+1$, one of the relations of $STVB_n$ which is not in $TVB_n$ and $SVB_n$. Then

   \begin{equation*}
f(\tau_i\gamma_j) = s_{e,\tau_i}s_{\overline{\tau_i},\gamma_j} = s_{e,\tau_i}s_{\rho_i,\gamma_j} = y_{i,i+1}\cdot\gamma_{j}
   \end{equation*}
     \begin{equation*}
f(\gamma_j\tau_i)= s_{e,\gamma_j}s_{\overline{\gamma_j},\tau_i}= s_{e,\gamma_j}s_{e,\tau_i}= \gamma_{j}\cdot y_{i,i+1}
    \end{equation*}
   $$r_{1,e}:=  y_{i,i+1}\cdot\gamma_{j}= \gamma_{j}\cdot y_{i,i+1}$$
  The remaining relations \( r_{1,\lambda} \), where \( \lambda \in \Lambda_n \), can be derived from this relation by conjugation with \( \lambda^{-1} \), yielding the same relation:
  $y_{ij}\gamma_k =\gamma_k y_{ij}$, by Lemma~\ref{gamma}. We have obtained (\ref{g5}).

Let us consider the last relation $r_2:=\tau_i=\rho_i\gamma_{i+1}\gamma_i\tau_i\gamma_i\gamma_{i+1}\rho_i$
\begin{align*}
f(\rho_i\gamma_{i+1}\gamma_i\tau_i\gamma_i\gamma_{i+1}\rho_i) & =s_{e,\rho_i}s_{\overline{\rho_i},\gamma_{i+1}}s_{\overline{\rho_i\gamma_{i+1}},\gamma_i}s_{\overline{\rho_i\gamma_{i+1}\gamma_i},\tau_i}\\
 & s_{\overline{\rho_i\gamma_{i+1}\gamma_i\tau_i},\gamma_i}s_{\overline{\rho_i\gamma_{i+1}\gamma_i\tau_i\gamma_i},\gamma_{i+1}}s_{\overline{\rho_i\gamma_{i+1}\gamma_i\tau_i\gamma_i\gamma_{i+1}},\rho_i}\\
 & =s_{e,\rho_i}s_{\rho_i,\gamma_{i+1}}s_{\rho_i,\gamma_{i}}s_{\rho_i, \tau_i}s_{e,\gamma_i}s_{e,\gamma_{i+1}}s_{e,\rho_i}\\
  &= \gamma_{i}\gamma_{i+1}(\rho_iy_{i,i+1}\rho_i)\gamma_i\gamma_{i+1}\\
  &= \gamma_{i}\gamma_{i+1}y_{i+1,i}\gamma_i\gamma_{i+1}
\end{align*}
and $f(\tau_i)=S_{e,\tau_i}=y_{i,i+1}$. Therefore, $r_{2,e}:=  y_{i,i+1}= \gamma_{i}\gamma_{i+1}y_{i+1,i}\gamma_i\gamma_{i+1}$.

Conjugating this relation by all representatives from $\Lambda_n$, we obtain (\ref{g6}).

Therefore, the monoid $STVP_n$ define by the relation (\ref{inverseexists})-(\ref{g6}).
\end{proof}
\subsection{The submonoid $STVH_n$ of monoid $STVB_n$}
\begin{figure}[ht]
  \centering
    \includegraphics[width=7cm]{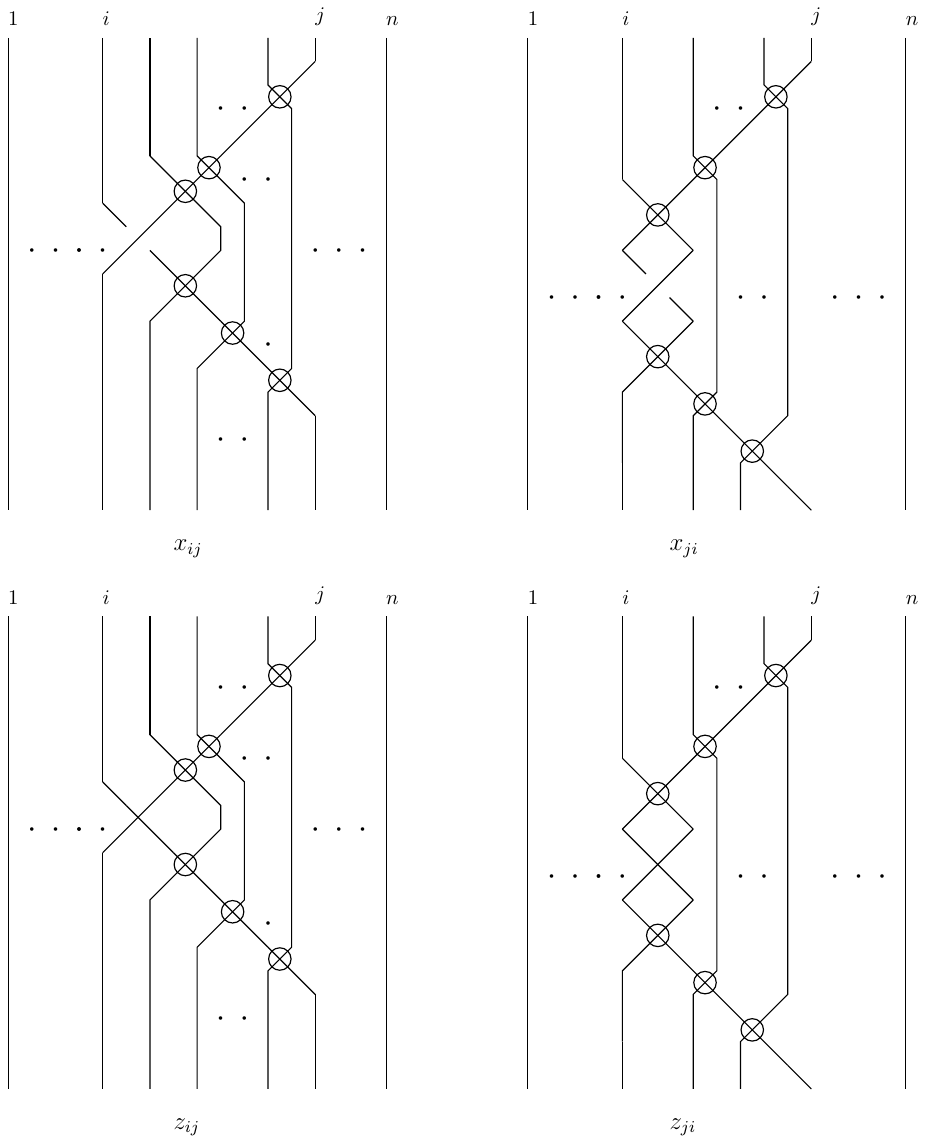}
        \caption{Generators of the monoid of $STVH_n$.}
        \label{genh}
        \end{figure}

Consider the following elements of  $STVB_n$:
$$x_{i,i+1}=\sigma_i,\text{   } x_{i+1,i}=\rho_i\sigma_i\rho_i, \text{   }  i=1,\ldots, n-1,$$
$$x_{i,j}=\rho_{j-1}\rho_{j-2}\cdots \rho_{i+1}\sigma_{i}\rho_{i+1}\cdots \rho_{j-2}\rho_{j-1},$$ 
$$ x_{j,i}=\rho_{j-1}\rho_{j-2}\cdots \rho_{i+1}\rho_i\sigma_i\rho_i\rho_{i+1}\cdots \rho_{j-2}\rho_{j-1},\text{   }  1 \leq i \leq j-1 \leq n-1,$$

$$z_{i,i+1}=\tau_i,\text{   } z_{i+1,i}=\rho_i\tau_i\rho_i, \text{   }  i=1,\ldots, n-1,$$
$$z_{i,j}=\rho_{j-1}\rho_{j-2}\cdots \rho_{i+1}\tau_{i}\rho_{i+1}\cdots \rho_{j-2}\rho_{j-1},$$ 
$$ z_{j,i}=\rho_{j-1}\rho_{j-2}\cdots \rho_{i+1}\rho_i\tau_i\rho_i\rho_{i+1}\cdots \rho_{j-2}\rho_{j-1},\text{   }  1 \leq i \leq j-1 \leq n-1.$$
Geometrically, these generators are shown in Figure~\ref{genh}.

\begin{theorem}\label{main2}
    The monoid $STVH_n$ admits a presentation with the generators $x_{kl}^{\pm 1}$, $z_{kl}$, $1 \leq k\neq l \leq n$, and
    $\gamma_j$, $1\leq j\leq n$. The defining relations are  (\ref{rel-inverse-b})-(\ref{rel-height-bb}) and as follows:
     \begin{align}
     x_{ij}x_{ij}^{-1} &= e\label{inverseexist2}\\
x_{ij}x_{kl} &=x_{kl}x_{ij}\label{2comm-clas}\\
x_{ik}x_{kj}x_{ik} &=x_{kj}x_{ik}x_{kj}\label{2classical2}\\
x_{ij}\gamma_k & =\gamma_kx_{ij}\label{g10}\\
x_{ij} &=\gamma_i\gamma_jx_{ji}\gamma_j\gamma_i\label{g11}\\
x_{ij}z_{kl} &=z_{kl}x_{ij}\label{2comm-clas3}\\
x_{ik}x_{kj}z_{ik} &=z_{kj}x_{ik}x_{kj}\label{2classical3}\\
z_{ik}x_{kj}x_{ik} &=x_{kj}x_{ik}z_{kj}\label{2classical4}\\
x_{ij}z_{ij} &= z_{ij}x_{ij}\label{2comm-clas4}\\
z_{ij}z_{kl} &=z_{kl}z_{ij}\label{2comm-clas2}\\
z_{ij}\gamma_k & =\gamma_k z_{ij}\label{g12}\\
z_{ij} &=\gamma_i\gamma_j z_{ji}\gamma_j\gamma_i\label{g13}
\end{align}
Where each distinct letter represents a unique index.
\end{theorem}
\begin{proof}
It is similar to the proof of Theorem~\ref{sTVP_n}.
\end{proof}
\subsection{The submonoid $M_n$ of monoid $STVB_n$}
\begin{theorem}\label{M_n}
    The monoid $M_n$ is defined by a presentation with the generators $\lambda_{kl}^{\pm 1}$, $z_{ij}$,~$1 \leq k\neq l \leq n$, and 
    $\gamma_j$, $1\leq j\leq n$. The defining relations are (\ref{rel-inverse-b})-(\ref{rel-height-bb}), (\ref{inverseexists})-(\ref{g4}), (\ref{2comm-clas2})-(\ref{g13}), and as follows:
     \begin{align}
\lambda_{ij}z_{kl} &=z_{kl}\lambda_{ij}\label{3comm-clas3}\\
\lambda_{ik}\lambda_{jk}z_{ji} &=z_{ji}\lambda_{ik}\lambda_{jk}\label{3classical3}\\
\lambda_{ik}\lambda_{ij}z_{kl} &=z_{kl}\lambda_{ik}\lambda_{ij}\label{3classical4}\\
\lambda_{ij}z_{ji} &= z_{ij}\lambda_{ij}\label{3comm-clas4}
\end{align}
Where each distinct letter represents a unique index.
\end{theorem}
\begin{proof}
It is similar to the proof of Theorem~\ref{sTVP_n}.
\end{proof}
\section{Some more results}

In this section, we establish additional results concerning the embedding of the monoid \( STVB_n \) into a group, as well as the extension of the \( TVB_n \) group representation to a representation of the \( STVB_n \) monoid.
\subsection{Embedding of monoid $STVB_n$ in a group $STVG_n$}
We demonstrate that the singular twisted virtual braid monoid on \(n\) strands embeds into a group, denoted by \(STVG_n\), which we refer to as the singular twisted virtual group on \(n\) strands. The group \(STVG_n\) contains a normal subgroup, \(STVPG_n\), consisting of singular twisted virtual pure braids.

In reference~\cite[Theorem 3]{E}, Keyman gave a technique employing elementary algebraic properties to demonstrate the embedding of specific types of monoids into groups. These monoids exhibit presentations akin to that of the singular braid monoids.

The free monoid on a set $A$ is denoted by $F^+(A)$. The Theorem is stated as follows:

\begin{theorem}~\cite{E}\label{ek}
    Let M be a monoid given by a presentation $[A \cup B \text{ } |\text{ } R]$
where $A = \{a_1,\ldots, a_n\}$, $B = \{b_1, \ldots, b_m\}$ and $R = R_1 \cup R_2 \cup R_3 \cup R_4$, where:
\begin{itemize}
    \item[(a)] \( R_1 \) comprises relations of the form $u = v$, where $u, v \in F^+(A)$,
\item[(b)] $R_2 = \{a_iu_i = u_ia_i = 1| \text{ for some } u_i \in F^+(A), \text{ for all } i = 1, \ldots, n\}$,
\item[(c)] $R_3$ consists of relations of the form $ub_j = b_ku$, for some $j, k = 1, \ldots, m$ and $u \in F^+(A)$,
\item[(d)] $R_4$ contains relations of the form $b_jb_k = b_kb_j$ , for some $j, k = 1, \ldots, m$.
\end{itemize}
Then M embeds in a group $G$ defined by the presentation $[A \cup B \cup \bar{B}\text{ } |\text{ } R \cup R']$, where $\bar{B}= \{\bar{b}_j \text{ } |\text{ } j =1, \ldots,m\}$and $R'$ consists of the following relations:
\begin{enumerate}

    \item $b_j\bar{b}_j = 1 =\bar{b}_jb_j$ ,
\item $u\bar{b}_j =\bar{b}_ku$ if $ub_j = b_ku \in R$, and
\item $\bar{b}_jb_k = b_k\bar{b}_j$ and $\bar{b}_j\bar{b}_k=\bar{b}_k\bar{b}_j$ if $b_jb_k = b_kb_j \in R$.
\end{enumerate}
\end{theorem}
In the presentation of a monoid  $M$, the invertible elements in set  $A$  can satisfy any relations among them. The elements in set  $B$  have no inverses (neither left nor right). If  $B = \emptyset$, then  M  is represented by the presentation $ [A \text{ } |\text{ } R_1 \cup R_2] $ and forms a group.

We demonstrate the embedding of monoid $STVB_n$ into a group utilizing the aforementioned theorem.
\begin{theorem}
    Let $n \in N$. The singular twisted virtual braid monoid $STVB_n$ embeds in the singular twisted virtual group $STVG_n$.
\end{theorem}
\begin{proof}
First we observe that $STVB_n$ satisfies the conditions in Theorem~\ref{ek}. Indeed, $STVB_n$ has a presentation in the form $[A \cup B \text{ } |\text{ } R_1 \cup R_2 \cup R_3 \cup R_4]$, where $A = \{ \sigma_i^{\pm 1}$, $\rho_i$ $(i=1, \dots, n-1)$, $\gamma_i$ $(i=1, \dots, n)\}$, $B = \{ \tau_i \text{ } |\text{ } 1 \leq i \leq n-1\}$, and
\begin{itemize}
    \item[(a)] $R_1$ comprises the defining relations: 
    
    $\sigma_i \sigma_j  = \sigma_j \sigma_i$,  $\rho_i \rho_{j} = \rho_{j} \rho_i$,  $\sigma_i \rho_{j}  = \rho_{j} \sigma_i$,  for $ |i-j| > 1$,

     $\gamma_j \rho_i  = \rho_i \gamma_j$,$ \sigma_i\gamma_j = \gamma_j\sigma_i$, for $j\neq i, i+1$,
     
    $\sigma_i \sigma_{i+1} \sigma_i  = \sigma_{i+1} \sigma_i \sigma_{i+1}$, $\rho_i \rho_{i+1} \rho_i  = \rho_{i+1} \rho_i \rho_{i+1} $,$ \rho_i \sigma_{i+1} \rho_i  = \rho_{i+1} \sigma_i \rho_{i+1}$,
   
 $ \gamma_{i+1} \rho_i  = \rho_{i} \gamma_i$, $ \rho_{i} \sigma_i \rho_{i}  = \gamma_{i+1} \gamma_i \sigma_{i} \gamma_i \gamma_{i+1}, \text{ for }   i=1,\ldots, n-1$,
 
 $\gamma_i \gamma_j  = \gamma_j \gamma_i \text{ for }   i,j=1,\ldots, n$.
     
     \item[(b)] $R_2$ contains the defining relations: 
     
$  \sigma_i \sigma_i^{-1}  = e$, $\rho_i^2=e$, $\gamma_i^2=e$, for $ 1 \leq i \leq n-1$.

 \item[(c)] $R_3$ includes the defining relations: 
 
 $\tau_i \rho_{j}  = \rho_{j} \tau_i  \text{ for }  |i-j| > 1,
\tau_i\gamma_j  = \gamma_j\tau_i  \text{ for }  j\neq i, i+1$,

$ \sigma_i \tau_{i}  = \tau_i \sigma_i, \rho_i \tau_{i+1} \rho_i  = \rho_{i+1} \tau_i \rho_{i+1}$,

$\sigma_i \sigma_{i+1} \tau_i  = \tau_{i+1} \sigma_i \sigma_{i+1},\sigma_{i+1} \sigma_{i} \tau_{i+1}  = \tau_{i} \sigma_{i+1} \sigma_{i}$

$\rho_{i} \tau_i \rho_{i}  = \gamma_{i+1} \gamma_i \tau_{i} \gamma_i \gamma_{i+1}  \text{ for }   i=1,\ldots, n-1.$

  \item[(d)] $R_4$ consists of the defining relations: 
  
     $\tau_i\tau_j = \tau_j\tau_i, \text{ for } |i-j| > 1.$
\end{itemize}
Therefore, the moniod $STVB_n$ embeds in the group $STVG_n$ as depicted below:

The singular twisted virtual group, \( STVG_n \), is defined as the group generated by the same generators as \( STVB_n \), together with \( \bar{\tau}_i \) for \( 1 \leq i \leq n-1 \), where \( \bar{\tau}_i \tau_i = e \) for all \( i \in \{ 1, \dots, n-1 \} \).
The defining relations consist of the following relations:
\begin{enumerate}

    \item The same monoid relations as $STVB_n$ with additional relations obtained by substituting $\bar\tau_i$ for $\tau_i$ in each $R_3$ relation, as shown below:
    
    $\bar{\tau}_i \rho_{j}  = \rho_{j} \bar{\tau}_i  \text{ for }  |i-j| > 1,
\bar{\tau}_i\gamma_j  = \gamma_j\bar{\tau}_i  \text{ for }  j\neq i, i+1$,

$ \sigma_i \bar{\tau}_{i}  = \bar{\tau}_i \sigma_i,  \rho_i \bar{\tau}_{i+1} \rho_i  = \rho_{i+1} \bar{\tau}_i \rho_{i+1}$,

$\sigma_i \sigma_{i+1} \bar{\tau}_i  = \bar{\tau}_{i+1} \sigma_i \sigma_{i+1},\sigma_{i+1} \sigma_{i} \bar{\tau}_{i+1}  = \bar{\tau}_{i} \sigma_{i+1} \sigma_{i}$

$\rho_{i} \bar{\tau_i} \rho_{i}  = \gamma_{i+1} \gamma_i \bar{\tau}_{i} \gamma_i \gamma_{i+1}  \text{ for }   i=1,\ldots, n-1.$

\item $\bar{\tau}_i\tau_i=e,  \text{ for }   i=1,\ldots, n-1.$

\item Lastly, relations obtained by substituting $\bar\tau_i$ for $\tau_i$ in each $R_4$ relation, as shown below:

$\bar{\tau}_i\bar{\tau}_j = \bar{\tau}_j\bar{\tau}_i, \bar{\tau}_i\tau_j = \tau_j\bar{\tau}_i \text{ for } |i-j| > 1.$
\end{enumerate}

\end{proof}

    We have established the group \( STVG_n \), which raises interest in the investigation of its subgroups. To achieve this, we define an epimorphism.

An epimorphism from \( STVG_n \) onto \( S_n \) can be defined as follows:

\begin{align*}
\pi \colon STVG_n \to S_n, & \quad \sigma_i \mapsto \rho_i, \quad \tau_i \mapsto \rho_i, \quad \rho_i \mapsto \rho_i, \quad i=1, 2, \ldots,n-1, \quad \\
& \gamma_j \mapsto e, \quad j=1, 2, \ldots, n.
\end{align*}

This mapping \(\pi\) forms an epimorphism, and the kernel of \(\pi\) denoted as $\ker(\pi)$, defines a subgroup. Specifically, the \(\ker(\pi)\) corresponds to the singular twisted virtual pure group \( STVPG_n \).

\begin{theorem}\label{sTVPG_n}
The group $STVPG_n$ has a presentation with the generators $\lambda_{kl},y_{kl}$, $1 \leq k\neq l \leq n$, and $\gamma_j$, $1\leq j\leq n$. The defining relations are  (\ref{rel-inverse-b})-(\ref{rel-height-bb}), and (\ref{comm-clas})-(\ref{g6}).
\end{theorem}
The proof proceeds in a manner analogous to that established for the submonoid \( STVP_n \).
\subsection{Extension of representation of $TVB_n$ to $STVB_n$}

V.G. Bardakov et al.~\cite{VNT} explore the extension of a given representation \( \varphi: B_n \to G_n \) of the braid group \( B_n \) into a group \( G_n \) to a representation \( \Phi: SM_n \to A_n \), where \( A_n \) is an associative algebra with \( G_n \) embedded in its group of units \( A_n^* \). The necessary conditions for this extension are stated in the following proposition.

\begin{proposition}~\cite{VNT}
Suppose \( \varphi: B_n \to G_n \) is a representation of the braid group \( B_n \), and let \( k \) be a field with elements \( a, b, c \in k \). Then the map \( \Phi_{a,b,c} : SM_n \to k[G_n] \), defined by the following action on generators:
\[
\Phi_{a,b,c}(\sigma_i^{\pm 1}) = \varphi(\sigma_i^{\pm 1}), \quad \Phi_{a,b,c}(\tau_i) = a\varphi(\sigma_i) + b\varphi(\sigma_i^{-1}) + ce, \quad i = 1, 2, \ldots, n - 1,
\]
establishes a representation of \( SM_n \) into \( k[G_n] \). Here, \( e \) represents the identity element in \( G_n \).
\end{proposition}

 Now the question is given a representation $\phi:TVB_n \to G_n$ of the twisted virtual braid group $TVB_n, n\geq 2$ into a group $G_n$. Can this representation be extended to a representation $\Phi: STVB_n \to A_n$, where $STVB_n$ is the singular twisted virtual braid monoid and $A_n$ is an associative algebra, in which the group of units contains $G_n$.

The answer to this question is affirmative. We can extend the representation from the group $TVB_n$ to monoid $STVB_n$. Suppose there is a representation $\phi:TVB_n \to G_n$ of the twisted virtual braid group $TVB_n, n\geq 2$ into a group $G_n$.

Then the map $\Phi: STVB_n \to \mathbb{K}[G_n]$, where $STVB_n$ is the singular twisted virtual braid monoid, $\mathbb{K}$ is a field and $ a, b, c \in \mathbb{K}$ which acts on the generators by the rule stated below. $$\Phi_{a,b,c}(\sigma_i)=\phi(\sigma_i), \Phi_{a,b,c}(\rho_i)=\phi(\rho_i),\Phi_{a,b,c}(\gamma_i)=\phi(\gamma_i),$$
$$\Phi_{a,b,c}(\tau_i)=a\phi(\sigma_i)+b\phi(\sigma_i^{-1})+ce.$$

We need to verify the combined relations and relations which involve only the generators \( \tau_i \). The relations (\ref{rel-height-vv2})-(\ref{rel-sss2}) have already been established in~\cite{VNT}.
So, now we need to check the relations~(\ref{rel-height-vv1})-(~\ref{1rel-twist-III}) only.
\begin{itemize}
    \item[(i)] $\tau_i \rho_{j} = \rho_{j} \tau_i$ for $|i-j| > 1$,

Apply, $\Phi_{a,b,c}$, we obtain the equality
\begin{align*}
    (a\phi(\sigma_i)+b\phi(\sigma_i^{-1})+ce) \phi(\rho_{j})  &= (a\phi(\sigma_i\rho_{j})+b\phi(\sigma_i^{-1}\rho_{j})+ce\phi(\rho_{j})\\
   &= (a\phi(\rho_{j}\sigma_i)+b\phi(\rho_{j}\sigma_i^{-1})+\phi(\rho_{j})ce \\
     &= \phi(\rho_{j})(a\phi(\sigma_i)+b\phi(\sigma_i^{-1})+ce)
\end{align*}
\item[(ii)] $\rho_i \tau_{i+1} \rho_i = \rho_{i+1} \tau_i \rho_{i+1}$  for $i=1,\ldots, n-2$,
Apply, $\Phi_{a,b,c}$, both the sides
\begin{align*}
   \phi(\rho_{i}) (a\phi(\sigma_{i+1})+b\phi(\sigma_{i+1}^{-1})+ce) \phi(\rho_{i})  &= (a\phi(\rho_{i}\sigma_{i+1}\rho_{i})+b\phi(\rho_{i}\sigma_{i+1}^{-1}\rho_{i})+ce\phi(\rho_{i}\rho_{i})\\
   &= (a\phi(\rho_{i+1}\sigma_i\rho_{i+1})+b\phi(\rho_{i+1}\sigma_i^{-1}\rho_{i+1})\\
   & \text{   }+\phi(\rho_{i+1}\rho_{i+1})ce \\
     &=\phi(\rho_{i+1}) (a\phi(\sigma_{i})+b\phi(\sigma_{i}^{-1})+ce) \phi(\rho_{i+1})
\end{align*}
\item[(iii)] $\tau_i\gamma_j = \gamma_j\tau_i$ for $ j\neq i, i+1$,
Apply, $\Phi_{a,b,c}$, we arrive at the equality
\begin{align*}
    (a\phi(\sigma_i)+b\phi(\sigma_i^{-1})+ce) \phi(\gamma_{j})  &= (a\phi(\sigma_i\gamma_{j})+b\phi(\sigma_i^{-1}\rho_{j})+ce\phi(\gamma_{j})\\
   &= (a\phi(\gamma_{j}\sigma_i)+b\phi(\gamma_{j}\sigma_i^{-1})+\phi(\gamma_{j})ce \\
     &= \phi(\gamma_{j})(a\phi(\sigma_i)+b\phi(\sigma_i^{-1})+ce)
\end{align*}
\item[(iv)]  $\rho_{i} \tau_i \rho_{i}  = \gamma_{i+1} \gamma_i \tau_{i} \gamma_i \gamma_{i+1}$ for $  i=1,\ldots, n-1$. Apply, $\Phi_{a,b,c}$, we get the equality
\begin{align*}
   \phi(\rho_{i}) (a\phi(\sigma_{i})+b\phi(\sigma_{i}^{-1})+ce) \phi(\rho_{i})  &= (a\phi(\rho_{i}\sigma_{i}\rho_{i})+b\phi(\rho_{i}\sigma_{i+1}^{-1}\rho_{i})+ce\phi(\rho_{i}\rho_{i})\\
   &= (a\phi(\gamma_{i+1} \gamma_i \sigma_{i} \gamma_i \gamma_{i+1})+b\phi(\gamma_{i+1} \gamma_i \sigma_i^{-1}\gamma_{i} \gamma_{i+1} )\\
   & \text{   }+\phi(\gamma_{i+1} \gamma_i \gamma_{i} \gamma_{i+1})ce \\
     &=\phi(\gamma_{i+1})\phi(\gamma_{i}) (a\phi(\sigma_{i})+b\phi(\sigma_{i}^{-1})+ce) \phi(\gamma_{i}) \phi(\gamma_{i+1})
\end{align*}

\end{itemize}

\section{Conclusion}
In conclusion, several intriguing problems remain open for further exploration. Firstly, can we determine the decomposition of the submonoids \( STVP_n \), \( STVH_n \), and \( M_n \)? Understanding their structure could provide deeper insights into their algebraic properties. Additionally, it is essential to explore whether the submonoids found are isomorphic or non-isomorphic, as this could significantly impact their classification and applications.
These questions offer promising avenues for further investigation and could lead to significant advancements in the study of algebraic structures related to twisted virtual braids and their extensions.
\section*{Acknowledgment}
We are thankful to Prof. V. G. Bardakov for his valuable discussions during the preparation of the paper on the twisted virtual braid group~\cite{VTKM}. His insights were particularly helpful in addressing issues that emerged concerning singular twisted virtual braids.
\section*{Funding}
The first author acknowledges the support of the Anusandhan National Research Foundation (ANRF), Government of India, through the ANRF Project CRG/2023/004921. The corresponding author expresses gratitude for the financial support from the University Grants Commission (UGC), India, in the form of a research fellowship under NTA Reference No. 191620008047. This work was also partially supported by the Department of Science and Technology (DST), Government of India, under the FIST program, Reference No. SR/FST/MS-I/2018/22(C).

\section*{Conflict of Interest}
The author declares no conflicts of interest.

\label{'ubl'}  
\end{document}